\documentclass[11pt,twoside]{article}
\usepackage{latexsym}
\usepackage{amssymb,amsbsy,amsmath,amsfonts,amssymb,amscd}
\usepackage{graphicx,color}
\usepackage{hyperref}
\newcommand{\commentout}[1]{}
\newcommand{\R}{\mathbb{R}}

\newcommand {\al} {\alpha}
\newcommand {\e}  {\varepsilon}

\newcommand {\Chi} {{\bf \raise 2pt \hbox{$\chi$}} }
\newcommand {\f}   {\frac}
\newcommand {\p}   {\partial}

\newcommand{\fer}{\eqref}
\newcommand{\dis}{\displaystyle}
\newcommand {\proof} {\noindent {\bf Proof}. }
\newcommand{\beq}{\begin{equation}}
\newcommand{\eeq}{\end{equation}}
\newcommand{\bal}{\begin{align}}
\newcommand{\bc}{\begin{cases}}
\newcommand{\ec}{\end{cases}}
\newcommand{\bea} {\begin{array}{rl}}
\newcommand{\eea} {\end{array}}
\newcommand{\bepa}{\left\{ \begin{array}{l}}
\newcommand{\eepa} {\end{array}\right.}
\definecolor{gre}{rgb}{0,0.7,0}

\newtheorem{theorem}{Theorem}[section]

\newtheorem{definition}[theorem]{Definition}

\newtheorem{prop}[theorem]{Proposition}

\newcommand{\qed}{{ \hfill
                     {\unskip\kern 6pt\penalty 500 \raise -2pt\hbox{\vrule\vbox to 6pt{\hrule width 6pt
                     \vfill\hrule}\vrule} \par}   }}
\title{\Large \bf Direct competition results from strong competition for limited resource}

\author{ Sepideh Mirrahimi\thanks{CMAP,  Ecole Polytechnique, CNRS, INRIA,  UMR CNRS 7641, F91128 Palaiseau cedex. Email: mirrahimi@cmap.polytechnique.fr}
\and
Beno\^it Perthame\thanks{
Universit\'e P. et M. Curie Paris 06, CNRS UMR 7598, Laboratoire Jacques-Louis Lions, F-75252 Paris cedex 05
}
\thanks{INRIA Paris-Rocquencourt, EPI BANG and Institut Universitaire de France. Email: benoit.perthame@upmc.fr}
\and Joe Yuichiro Wakano\thanks{
Meiji Institute for Advanced Study of Mathematical Sciences and PRESTO, Japan Science and Technology Agency. Email: joe@meiji.ac.jp}}

\date{\today}

\begin{document}
\maketitle
\pagestyle{plain}
\pagenumbering{arabic}

\begin{abstract}
We study a model of competition for resource through a chemostat-type model where species consume the common resource that is constantly supplied. We assume that the species and resources are characterized by a continuous trait. As already proved, this model, although more complicated than the usual Lotka-Volterra direct competition model, describes competitive interactions leading to concentrated distributions of species in continuous trait space.

Here we assume a very fast dynamics for the supply of the resource and a fast dynamics for death and uptake rates.   In this regime we show that factors that are independent of the resource competition become as important as the competition efficiency 
and that the direct competition model is a good approximation of the chemostat. 
Assuming these two timescales allows us to establish a mathematically rigorous proof showing that our resource-competition model with continuous traits converges to a direct competition model. We also show that the two timescales assumption is required to mathematically justify the corresponding classic result on a model consisting of only finite number of species and resources (MacArthur, R. Theor. Popul. Biol. 1970:1, 1-11).

This is performed through asymptotic analysis, introducing different scales for the resource renewal rate and the uptake rate. The mathematical difficulty relies in a possible initial layer for the resource dynamics. 

The chemostat model comes with a global convex Lyapunov functional. We show that the particular form of the competition kernel derived from the uptake kernel, satisfies a positivity property which is known to be necessary for the direct competition model to enjoy the related Lyapunov functional.
\end{abstract}


\section{Introduction}
\label{sec:intro}

In order to faithfully capture the dynamics of species competition, a detailed mechanism such as species-resource dynamics should be explicitly modeled.  On the other hand, specifying a particular mechanism narrows the applicability of a model.  This is why conceptual and simple Lotka-Volterra type models of direct competition are widely and intensively studied.  Then a natural question arises: Can we reduce a species-resource dynamics model  into a Lotka-Volterra model?  If so, how are parameters (e.g., intrinsic growth rates, competition coefficients) determined by those of the original mechanistic model?  
\\

These questions are central when studying the evolutionary origin of species diversity which has been one of the most important problems in evolutionary ecology. Individuals with similar traits (e.g. body size or shape) face strong competition, which
might lead to adaptive radiation. Using models of competition-driven speciation, several theoretical
studies have shown that the species distribution in continuous trait space often evolve
toward single peak or multiple peaks that are distinct from each other  (see \cite{AS.SE:95, OD.PJ.SM.BP:05, LD.PJ.SM.GR:08, GB.BP:08, PJ.GR:09,SM.BP.JW:10}).  Although many previous results are derived from a simple model that assumes direct species competition of Lotka-Volterra type, competitive interaction among species generally occurs in competition for resource such as prey or nutrient.  For example, birds with similar beak shapes are in competition because they utilize similar food resource. Another example is found in ecological stoichiometry where consumer species with similar C:N:P (carbon: nitrogen: phosphorus ratio)
requirements experience competition for nutrient supply with their optimal C:N:P ratio (see \cite{RS.JE:02}).  Thus, the competitive interaction should be modeled not directly but implicitly through the resource. Recently, we have proposed such a mechanistic model of species competition where species (with continuous trait) compete for common resource (which is characterized by another continuous 'trait' or 'type') which is constantly supplied (chemostat-type model)  (see \cite{SM.BP.JW:10}).  We have shown that evolutionary branching, which is parallel to adaptive radiation, is possible in this chemostat-type model as well, supporting the results previously obtained by direct competition models.
\\

This paper aims the reduction of a mechanistic model of species-resource dynamics to a more conceptual model of direct species competition.  In Section  \ref{sec:model}, we describe our chemostat model of species-resource dynamics and a brief summary of our reduction result.  Section \ref{sec:ddc} is devoted to a rigorous description of this result and the proof.  In section \ref{mutations}, we extend our result to include a system with mutation.  In section \ref{sec:lyapunov}, we study the Lyapunov functional property which is useful for evolutionary stability analysis.  In section \ref{sec:numerics}, we show some numerical examples to study how the reduced dynamics approximates the original dynamics.  In section \ref{sec:macarthur}, we show that the classic result on a model of discrete number of species \cite{macarthur} was not fully mathematical and discuss how our framework can be used to solve this problem.

\section{Model and a brief summary of the results}
\label{sec:model}

In general, living organisms require several different resource to survive and reproduce.   They might undergo strong competition for some resource, while some other resource might be sufficiently supplied so that it might not lead to competition at all.  In case of phytoplankton species, resources are light and nutrients.  Although tremendous amount of light is supplied from the sun, the competition for light still exists.  Nutrient supplied from the bottom is also a limiting resource, and this dynamics might have much slower timescale compared to the competition for light (but slower dynamics do not imply less importance).  There might exist some features that determine growth or death rates independently of resource competition. Here we consider a model that allows continuous distribution of species trait and resource type (by an abuse we sometimes refer to it as a trait too)  for two reasons: first, a continuous model is more compatible with standard evolutionary models of physiological traits mentioned earlier; second, the mathematical analysis of a continuous model seems more general than a model with the finite and given number of species and resources (Note that our mathematical treatment allows distribution to converge to Dirac masses.  Thus, although the rigorous proof is not given here, it implies that our model includes relatively simple ODE models as special cases.)

Inspired by such phytoplankton ecology, our mechanistic model for a chemostat with continuous resource supply and continuous consumer population trait is given by
\begin{equation}
\label{basic}
\begin{cases}
\f{\p}{\p t}n(x,t)={n(x,t)}\big[b(x)-d(x) + \int K(x,y)R (y,t) dy\big],
\\[4mm]
 \f{\p}{\p t} R (y,t)= {m(y)} \big[R_{\rm in}(y) - R (y,t) \big]- R (y,t)\int K(x,y)n (x,t)dx,
\end{cases}
\end{equation}
where  $n(x,t)\geq 0$ is the consumer species density with trait $x \in \R$ and $R(y,t)\geq 0$ is the common resource density with trait $y \in \R$.  We assume that a death rate $d(x)$ depends on species trait $x$.  Reproduction rate consists of two parts: $b(x)$ denotes a basic reproduction rate independent of the competition for the common resource, and $\int K(x,y)R (y,t) dy$ denotes a reproduction rate coming from consuming the common resource $R$.  The resource $R$ is supplied into the system with a rate $m(y)$ so that $R(y,t)$ approaches $R_{\rm in}(y)>0$. These drive the system (a chemostat type model).  Consumption for resource is denoted by a trait dependent uptake kernel $K(x,y) \geq 0$ that defines how species with trait $x$ depends on resource with trait $y$.  All coefficients ($b(x),d(x),K(x,y),m(y),R_{\rm in}(y)$; although they are functions of traits, we call them 'coefficients') can play an important role in the full system. 

 The regime when direct competition is close to \eqref{basic} arises in some limiting situation.  We assume that the dynamics has the following three different timescales (notice however that time can always be rescaled to fix one of them, and as we explain in section \ref{sec:macarthur}, only two timescales cannot do the work)
\\[4pt]
$\bullet$ {\em Very fast dynamics for the common resource supply.}  This imposes that, in the absence of consumer species, resource distribution $R$ approaches $R_{\rm in}$.
\\
$\bullet$ {\em  Fast dynamics for resource-consumption and the resulting reproduction and death.}  In other words, common resource dependent coefficients  are large. This implies that competition for the common resource is a primary factor for the growth of consumer species population.
\\
$\bullet$ {\em Slow dynamics.}	Population dynamics of $n$ is also affected by a basic reproduction rate $b(x)$ and a death rate $d_{\rm slow}(x)$, both of which are independent of the common resource.
\\

In order to include these assumptions in the equations, we scale accordingly the coefficients
\\[4pt]
$\bullet$ $m$ is of order $\e^{-2}$,
\\
$\bullet$ $K$ is of order $\e^{-1}$,
\\
$\bullet$ $d(x)=    \int K(x,y)R_{\rm in}(y) dy + d_{\rm slow}(x)$,
\\
and we define a neat growth/death rate (that we take independent of $\e$ for simplicity)
\\
$\bullet$ $a(x) =   b(x) - d_{\rm slow}(x)$.
\\
\\
In other words, direct competition is a good approximation of \eqref{basic}  when the death rate is suited to counterbalance fast growth based on resource-consumption (otherwise, the population will grow infinitely).  For this purpose, in the third bullet, we divide the death rate into two parts.  The first term represents fast dynamics for death events that compensates fast resource-consumption.  The second term, $d_{\rm slow}(x)$, represents the remaining slow dynamics for death events independent of the competition for resource.  Note that this decomposition of $d(x)$ means that, while the common resource $R$ is at its equilibrium $R_{\rm in}$, the dynamics of the consumers are determined only by the parameters which are independent of the common resource, that is why we introduce the neat rate  $a(x)$.


These result in a new formulation of \eqref{basic} that is written  as
\begin{equation}
\label{full}
\begin{cases}
\f{\p}{\p t} n_\e(x,t)={n_\e(x,t)}\big[a(x) + \f 1 \e \int K(x,y)\big(R_\e (y,t)- R_{\rm in}(y) \big) dy\big],
\\[4mm]
\f{\p}{\p t}  R_\e (y,t)= \f {m(y)}{\e^2} \big[R_{\rm in}(y) - R_\e (y,t) \big]- \f 1 \e R_\e (y,t)\int K(x,y)n_\e (x,t)dx,
\end{cases}
\end{equation}
where new coefficients satisfy $m=O(1)$, $K=O(1)$, and $a=O(1)$.

We complement the system with initial data 
$$
 n_\e(x,0)=n^0(x) \geq 0 ,\qquad R_\e(y,0)=R^0(y) \geq 0, 
$$

As the very fast and fast dynamics become infinitely fast, $R$ should
approach $R_{\rm in}$, but it is not clear which distribution $n$ should approach.  When $R$ is not close to $R_{\rm in}$, the reproduction based on resource-consumption is dominant (fast dynamics with order $\e^{-1}$) so the dynamics of $n$ is primarily determined by $R-R_{\rm in}$.  However, as $R-R_{\rm in}$ approaches to zero, resource-dependent reproduction might become no longer dominant.

In this paper, we will show that an approximation of this dynamics is given by the direct competition model
\beq
\f{\p}{\p t} n(x,t)={n(x,t)}\big[a(x) - \int c(x,x') n(x',t) dx' \big] ,
\label{comp_model}
\eeq
where the binary competition kernel $c(x,x')$ is related to the coefficients in \eqref{full}.  We will show that it is given by
$$
c(x, x') =  \int K(x,y)  \f{R_{\rm in}(y)}{m(y)}  K(x',y) dy.
$$
This equation can be seen as a conversion formula of kernels.  The reduced competition-kernel $c(x,x')$ is very sensitive to $R_{\rm in}(y)/m(y)$ as well as to resource-consumption kernel $K(x',y)$, but not to $R_{\rm in}(y)$ or $m(y)$ separately.  
Note that this kernel is of order one (recall that resource-consumption kernel is of order $\e^{-1}$ before rescaling) so the reduced direct competition model describes the slow dynamics of $n(x,t)$.  We will later show that the inclusion of mutation in consumer species' trait does not essentially change the result.

One biologically interesting point of the equation \eqref{comp_model} is that the growth rate of consumer species in the absence of competition for resource, which is described by $a(x)$, still has considerable effect on the long-term dynamics of consumer distribution $n$.  
This is because as $R$ becomes very close to $R_{\rm in}$,  the dominant factor
(fast dynamics) becomes weaker and a basic growth rate $a(x)$ becomes as important as the competition for resource.  In this limit, we recover a Lotka-Volterra type competition dynamics from a mechanistic model of a chemostat type.  It is natural that we do not have the reduced dynamics of the resource, $R$, because the reduction is possible when the detailed dynamics of the resource can be neglected and when the effects of competition for the resource can be put into the new competition coefficients in a Lotka-Volterra type system  (i.e. $c(x,x')$).

As already raised in \cite{PA.CR.RD:08}, another practical question of our interest is the choice of the competition kernel; usually a Gaussian choice is favored, based on heuristic arguments which promotes the coexistence of infinitely many species (continuous distribution over continuous trait axis). Not only this has a particular effect as shown by the structural instability of continuous distribution in  result \cite{MG.GM:05}, but also our result shows that this choice is arbitrary because the conversion of resource-consumption kernel into direct-competition kernel does not guarantee that a Gaussian function is converted into another Gaussian function. Take 
$$
K(x,y) = e^{-\alpha (x-y)^2} , \qquad R_{\rm in} (y) = e^{- \beta y^2}, \qquad m(y) \equiv 1.
$$
Then, choosing $\gamma $ such that $\gamma (2\alpha+\beta) = \alpha^2$ so as to complete a square in the expression
$$\bea
c(x,x')&=  e^{-[\al x^2 + \al (x')^2 -\gamma ({x+x'})^2]}  \int e^{-[(2\al+\beta) y^2 - 2 \al (x+x')y + \gamma ({x+x'})^2]}  dy
\\
&= e^{-[\al x^2 + \al (x')^2 -\gamma ({x+x'})^2]} \; \sqrt{\f{ \pi}{2\al + \beta}}, 
\eea$$ 
we find that $c(x,x')$ is not gaussian because it is proportional to 
$$
e^{-[\al x^2 + \al (x')^2 -\gamma ({x+x'})^2]} \neq e^{- (\al - \gamma)  (x-x')^2} \quad \text{ when } \beta \neq 0.
$$

The present study is not the first attempt to derive direct competition Lotka-Volterra equations from a model of  species competing for resources. MacArthur \cite{macarthur} proposed a simpler approach where he also chooses a Lotka-Volterra equation for the resources; he chooses a discrete setting rather than continuous but this is irrelevant for our purpose. With our notations MacArthur's model is described as a prey-predator model
\beq
\begin{cases}
\f{\p}{\p t}n(x,t)=b(x) n(x,t)\big[ \int K(x,y)R (y,t) dy- d(x) \big],
\\[4mm]
 \f{\p}{\p t} R (y,t)= R(y,t) \big[m(y) (1 - \f{R (y,t) }{R_{cc}(y)}) - \int K(x,y)n (x,t)dx \big],
\end{cases}
\label{macarthur}
\eeq
Assuming $b$ is of order $\e$ and rescaling time accordingly, we end up with a rescaled equation that is close to \eqref{full}. This seems to provide a simpler approach to obtain the aymptotic behavior. However, the shortcoming is that the limiting equation is NOT under the form \eqref{comp_model} but more complicated and with different nonlinearities (e.g., see \eqref{ma_comp} in section \ref{sec:macarthur}.) Thus, more rigorous analysis is necessary.
\\

Mathematical notations and results will follow in the next three sections.  Section \ref{sec:ddc} is devoted to a rigorous derivation of our convergence result. It goes through asymptotic analysis and compactness estimates that allow us to show that $n_\e$ converges to $n$ and thus $n$  is the population density. In other words our theory is not a perturbation analysis on $n_\e$ while it is on $R_\e \approx R_{\rm in}$. Mutations are not included in the model equations \eqref{basic}--\eqref{comp_model}. In section \ref{mutations}, we show how to include them and which new rescaling this leads to perform.  Another remarkable property of the chemostat is the Lyapunov functional property already used in \cite{SM.BP.JW:10}. The asymptotic form of this Lyapunov functional is related to that for the direct competition model in section \ref{sec:lyapunov} and in particular we explain why, in this context,  the corresponding operator is always nonnegative.

\section{Derivation of the direct competition}
\label{sec:ddc}

In order to state our result we need some assumptions
\beq
n^0, \; R^0, \; R_{\rm in} \in L^1\cap L^\infty(\R), \qquad\int \big|\f{\p}{\p x} n^0\big|\,dx\leq C_1,
\label{as:ID}
\eeq
\beq
a_M:= \| a \|_\infty< \infty,  \qquad {\underline m} = \min_{y} m(y) >0, \qquad {\overline m} = \max_{y} m(y) < \infty,
\label{as:coef0}
\eeq
\beq
K\in BUC(\R \times \R), \qquad  K_M := \sup_{x,y} K(x,y) < \infty,
\label{as:coef1}
\eeq
\beq
\big|\f{\p}{\p x}K(x,y)\big|\leq C_2,\qquad \big|\f{\p}{\p x} a(x)\big|\leq C_3.
\label{as:coef2}
\eeq
where $BUC$ stands for the linear space of Bounded Uniformly Continuous functions.

With these in mind, we  establish that the limiting problem for \eqref{full} is the direct competition model
\begin{prop} As $\e$  tends to $0$ in \eqref{full}, $R_\e$ converges strongly to $R_{\rm in}$ in $L^1 \big((0,T)\times \R \big)$ and $n_\e$ converges strongly in $C \big((0,T); L^1(\R)  \big)$ to $n$, the solution to \eqref{comp_model} with the initial data $n^0$ and the competition kernel
\beq
c(x, x') =  \int K(x,y)  \f{R_{\rm in}(y)}{m(y)}  K(x',y) dy.
\label{comp_kernel}
\eeq
\label{prop:ntoc}
\end{prop}

It is more technical to improve the space-time convergence on $R_\e$. Indeed it leaves place for an initial layer which is needed because, in general,   its limit $R_{\rm in}$ is incompatible with the initial data $R^0_\e$. 
\\

\proof To begin with, we re-write the equation on $R_\e$  in \eqref{full} as
$$
R_\e (y,t) - R_{\rm in}(y) = - \e  \f{R_\e(y,t)}{m(y)} \int K(x,y)n_\e (x,t)dx \,- \,\f{\e^2}{m(y)}  \f{\p}{\p t} R_\e (y,t). 
$$ 
We inject this to compute the growth rate for $n_\e$ 
$$
 \f 1 \e \int K(x,y)\big(R_\e (y,t) - R_{\rm in} \big) dy= -  \int K(x,y) \f{R_{\rm in }(y,t)}{m(y)}  \int K(x',y) n_\e(x',t) dx' dy+I_\e(x,t)
$$
$$
\qquad \quad \qquad = - \int c(x,x') n_\e(x',t) dx' +I_\e(x,t),
$$
as we use $c(x, x') $ as  defined in \eqref{comp_kernel} and with 
$$
I_\e (x,t) = - \e \f{\p}{\p t} \int \f{K(x,y) }{m(y)} R_\e (y,t) dy+\int K(x,y) \f{R_{\rm in}(y)-R_{\e}(y,t)}{m(y)}  \int K(x',y) n_\e(x',t) dx'dy.
$$

Therefore the equation on $n_\e$ can also be written
\beq
\f{\p}{\p t} n_\e(x,t)={n_\e(x,t)}\Big[a(x) - \int c(x,x') n_\e(x',t) dx' +I_\e(x,t) \Big].
\label{ddc:reqn}
\eeq
Throwing away the terms in $\e$, we find formally \eqref{comp_model}.
\\

To justify rigorously this limit, we need (i) estimates for compactness, (ii) explain how we can pass to the limit. 
\\[4pt]
\noindent  {\em (i) Estimates.} We define 
$$
M_\e(t) : = \int n_\e(x,t) dx, \qquad N_\e(t) := \int \big|R_\e(y,t) -R_{\rm in}(y) \big| dy, 
$$
We then notice from  the equation on $R_\e$ in \fer{full} that we have 
$$R_\e(y,t)\leq R_{\rm in}(y),\qquad \text{for all $y$ such that $R^0(y)\leq R_{\rm in}(y)$,}$$
and
$$R_\e(y,t)\leq R_{\rm in}(y)+\left(R^0(y)-R_{\rm in}(y)\right)\exp(-\f{\underline m t}{\e^2}),\qquad \text{for all $y$ such that $R^0(y)> R_{\rm in}(y)$.}$$
It follows in particular that $N_\e$ and $\|R_\e(\cdot,t)\|_{L^1}$ are uniformly bounded in $\e$ by a constant denoted by $R_M$.
Next we write from the equation on $n_\e$ in \fer{full}
$$
\f{d}{dt}M_\e\leq  M_\e \left[ a_M +\f{K_M}{\e}\exp\left(-\f{\underline m t}{\e^2}\right)\int \left(R^0(y)-R_{\rm in}(y)\right)_+dy \right].
$$
We deduce that, for $\overline C$ some constant independent of $\e$,
$$M_\e(t)\leq M_\e(0) \exp\left(a_M t+\overline C\e\right),$$
and thus $M_\e$ is uniformly bounded in $\e\leq 1$ on each time interval $[0,T]$ by a constant denoted by $M_T$.
Next, we define and write from the equation on $R_\e$ in \fer{full}
$$
\f{d}{dt} N_\e(t)  + \f {\underline m} {\e^2}N_\e(t)   \leq  \f {R_M} {\e}  K_M M_\e(t) .
$$ 
It follows from this differential inequalitiy and the fact that $M_\e(t)$ and $N_\e(t)$ are uniformly bounded that $\f 1 \e \int_0^t N_\e(s)ds$ is also uniformly bounded in $\e\leq 1$ on each time interval $[0,T]$. One step further, from the bound on $M_\e$ and from the equation on $N_\e$, one also concludes that  
\beq\label{Nee}
\f{N_\e(t)}{\e} \leq \f{N(0)}{\e} e^{- {\underline m} t/\e^2}+ \f {R_M K_M M_T} {\underline m}  , \qquad \forall t \in [0,T].
\eeq
\\
\noindent  {\em (ii) Compactness in $t$.}
From this, and going back to the equation on $n_\e$, we conclude that, still for $t \in [0,T]$
\beq
 \int \left|\f{\p n_\e(x,t)}{\p t} \right| dx \leq a_M C(T) + K_M \f{N(0) }{\e} e^{- {\underline m} t/\e^2}+ K_M C(T).
\label{ddc:est1}
\eeq
This Lipschitz estimate in time gives us that $n_\e$ is uniformly equicontinuous in $C\big((0,T); L^1(\R)\big)$.
\\
\noindent {\em (iii) Compactness in $x$.} Dividing equation \eqref{full} by $n_\e$ and differentiating in $x$, we find
$$
\f{\p^2}{\p x \p t}\ln n_\e(x,t)=\f{\p}{\p x}\, a(x)+\f{1}{\e}\int \f{\p}{\p x} \,K(x,y)\, \left(R_\e(y,t)-R_{\rm in}(y)\right)dy.
$$
From \eqref{as:coef2} and \eqref{Nee} we obtain that, for some constant $D_1$
$$
\big|\f{\p^2}{\p x \p t}\ln n_\e(x,t)\big|\leq \,D_1 \,+\,   C_2 \f{N(0)}{\e} e^{- {\underline m} t/\e^2}.
$$
It follows that,  for $\e\leq 1$ and for some constant $D_2$
$$
\big|\f{\p}{\p x }\ln n_\e(x,t)\big|\leq \big|\f{\p}{\p x }\ln n^0(x)\big|+ D_2 t,
$$
$$
\int \big|\f{\p}{\p x } n_\e(x,t)\big|dx\leq \int \f{\big|\f{\p}{\p x } n^0(x)\big|}{n^0(x)}\,n_\e(x,t)\,dx+ D_2 t \int n_\e(x,t) dx. 
$$
Moreover, from the first equation of \eqref{full} and \eqref{Nee} we obtain that, for some constant $D_3$,
$$
\big|\f{\p}{\p t} \ln n_\e(x,t)\big|\leq D_3\,+\,   K_M \f{N(0)}{\e} e^{- {\underline m} t/\e^2},
$$
and thus: $\ln n_\e(x,t)\leq  \ln n^0(x) +D_4 \,t$, for $\e\leq 1$ and some constant $D_4$, which gives
$$
 \f{n_\e(x,t)}{n^0(x)}\leq \exp(D_4 t).
$$
Using the above arguments and \eqref{as:ID}, we conclude that, for some constant $D(T)$ and $0\leq t \leq T$,
$$\int \big|\f{\p}{\p x } n_\e(x,t)\big|\,dx\leq \exp(D_4 t)\int \big|\f{\p}{\p x } n^0(x)\big|\,dx+D_2 t \int n_\e(x,t) dx\leq C_1\exp(D_4 t)+D_2 tM_\e(t)\leq D(T). 
$$
\\

Using this and step (ii) we conclude from the Arzela-Ascoli Theorem that,  after extracting a subsequence, $n_\e$ converges in $C\big((0,T); L^1(\R)\big)$.
\\
\\
{\em (iv) Passing to the limit.} It remains to pass to the weak (distribution) limit in the equation on $n_\e$ written as \eqref{ddc:reqn}.  Because of the strong convergence of $n_\e$ and \eqref{Nee}, the only difficulty is to pass to the limit in the term
\begin{align}\nonumber
\e n_\e(x,t) \f{\p}{\p t} \int \f{K(x,y) }{m(y)} R_\e (y,t) dy&= \e  \f{\p}{\p t} \big[ n_\e(x,t) \int \f{K(x,y) }{m(y)} R_\e (y,t) dy \big]\\\nonumber
& - \e \int \f{K(x,y) }{m(y)} R_\e (y,t) dy \f{\p}{\p t} n_\e(x,t).
\end{align}
The first term converges weakly to $0$ (multiplying by a test function, after integration by parts all the terms are bounded and multiplied by $\e$). For the second term, we use the estimate \eqref{ddc:est1} to reach the conclusion that it converges in $L^1$ to $0$, thanks again to the multiplying factor $\e$.

Since uniqueness of weak solutions for the direct competition model is proved in \cite{LD.PJ.SM.GR:08}, we conclude that the full family $n_\e$ converges and not only subsequences.
\qed

\section{Mutations}
\label{mutations}

In the present context of population equations structured by a physiological trait, there are several possibilities to represent mutations that have been widely used \cite{MR2041498, OD.PJ.SM.BP:05, LD.PJ.SM.GR:08}. Integral operators or diffusion operators can be derived from stochastic individual based models; see for instance \cite{NC.RF.SM:06, NC..RF.SM:08bis}. Mathematically they have many similar properties, in particular asymptotic analysis can be carried out using similar methods in the regime of small mutations for a long time of observation leading to speciation phenomena \cite{ GB.BP:07,  GB.SM.BP:09, AL.SM.BP:10} although stochastic individual methods are also used for the same purpose \cite{MR2250806,NC.SM:11,SM.lisbon}. 

Here, we restrict ourselves to representing the mutations by a diffusion operator with intensity $\mu^2$, a very small rate. According to the mathematical theory in \cite{OD.PJ.SM.BP:05, GB.BP:08}, it is natural to rescale time so that the new timescale unit is $1/\mu$, very long timescale compared to the original 'generation' timescale. This leads to the competition model with mutations
\begin{equation}
\label{full_mut}
\begin{cases}
\mu \f{\p}{\p t} n_\e(x,t)={n_\e(x,t)}\big[a(x) + \f 1 \e \int K(x,y)\big(R_\e (y,t)- R_{\rm in} \big) dy\big] + \mu^2 \Delta n_\e,
\\[4mm]
\mu \f{\p}{\p t} R_\e (y,t)= \f {m(y)}{\e^2} \big[R_{\rm in}(y) - R_\e (y,t) \big]- \f 1 \e R_\e (y,t)\int K(x,y)n_\e (x,t)dx,
\end{cases}
\end{equation}

Because both the uptake rate $\e$ and the mutation rate $\mu$ can be considered as small, it is interesting to see how the solutions behave in the different regimes of smallness.

We can  formally follow  the lines of the analysis in section \ref{sec:ddc} and find direct competition with mutations as the limit of infinitely fast uptake rate 
\beq
\mu \f{\p}{\p t} n(x,t)={n(x,t)}\big[a(x) - \int c(x,x') n(x',t) dx' \big] + \mu^2 \Delta n .
\label{comp_mut}
\eeq

One can also consider the small mutation rate and, assuming the initial data is well prepared, use the Hopf-Cole unknowns
$$
u_\mu = \mu \ln(n_\e), \qquad  u^0_\mu = \mu \ln(n^0).
$$
This leads now to study the equation on the limit $u_\e(x,t)$ of $u_\mu(x,t)$ and on the limit measure $n_\e$ as $\mu \to 0$. Following formally \cite{OD.PJ.SM.BP:05,GB.BP:08}, one finds the constrained Hamilton-Jacobi equation
$$\begin{cases}
 \f{\p}{\p t} u_\e = a(x) + \f 1 \e \int K(x,y)\big(R_\e (y,t)- R_{\rm in} \big) dy + | \f{\p u_\e}{\p x} |^2, \qquad \dis \max_{x\in \R} u_\e(x,t)=0 \quad \forall t \geq 0, 
\\[4mm]
 \f {m(y)}{\e^2} \big[R_{\rm in}(y) - R_\e (y,t) \big]- \f 1 \e R_\e (y,t)\int K(x,y)n_\e (x,t) dx =0 ,
\end{cases}$$
and $n_\e$ is supported in the zeroes of $u_\e$ (usually points).

This can be further analyzed in the regime $\e$ small and leads to 
$$
 \f{\p}{\p t} u = a(x) +  \int c(x,x') n(x',t) dx' + \big| \f{\p u }{\p x}\big|^2, \qquad \max_{x\in \R} u(x,t)=0 \quad \forall t \geq 0, 
$$
still with the information that the measure $n(t)$ is supported by the points where $u(t)$ vanishes. This is also the limit of \eqref{comp_mut} in the regime of small $\mu$. 
\\

We can conclude that the formal analysis of small mutations, both for the chemostat model and the direct competition model, is compatible with the multiscale analysis in section \ref{sec:ddc} that relates both models. Of course a rigorous analysis again requires mathematical developments,  in particular appropriate estimates,which are beyond the scope of the present paper. This is the regime of interest that we choose later in section   \ref{sec:numerics} for numerical illustration.  In this regime, we recover also, at a population level,  the results of standard adaptive dynamics \cite{SG.EK.GM.JM:98, OD:04} and the references therein.

\section{Lyapunov functionals for the two competition models}
\label{sec:lyapunov}

The large time behaviour of solutions to \eqref{full} has been studied in  \cite{SM.BP.JW:10}. It is proved that steady states cannot be globally positive  $(\overline n_\e, \overline R_\e)$ and they have to be either Dirac masses (speciation) or continuous with a support of $n_\e$ that is small enough (except when the solution goes extinct that might happen if renewal/death  rates are too  high). When continuous, the steady state globally attracts all trajectories if they satisfy a particular sign property characterized the so-called {\em Evolutionary Stable Distribution} (ESD in short),  a notion that we recall below. The proof relies on a Lyapunov functional which has also been proved to exist for the direct competition model in \cite{PJ.GR:09} under the condition of positivity of a certain operator (see also \cite{NC.PJ.GR:11} for a similar study of a model with finite number of resources and consumer traits).

Our interest here is to understand the relation between the two Lyapunov functionals,  for \eqref{full} and \eqref{comp_model}, and to understand why the operator positivity condition is important for direct competition models but does not appear in models of competition for resources.

We first recall that  an ESD is defined as
\begin{definition} [Evolutionary Stable Distribution for \eqref{comp_model}, \cite{PJ.GR:09}]  A nonnegative bounded measure $\overline n$ is called an ESD for the direct competition equation \eqref{comp_model} 
\beq \label{steady_dc1}
a(x)-   \int c(x,x') \overline n(x')dx'= 0 \qquad \forall x\in {\rm Supp}\; \overline n, 
\eeq
\beq
\label{steady_dc2}
\quad a(x) - \int c(x,x')\overline n(x')dx' \leq 0 \qquad \forall x\in \mathbb{R} \setminus {\rm Supp}\; \overline n.
\eeq 
\end{definition}
\begin{definition} [Evolutionary Stable Distribution for \eqref{full}, \cite{PJ.GR:09, SM.BP.JW:10}]  A state  is $(\overline n_\e, \overline R_\e)$  is called an ESD for the competition for resources equation \eqref{full} if $\overline n_\e$ is a nonnegative bounded measure and
\beq \label{steady1}
a(x)+ \f 1 \e \int K(x,y) [\overline R_\e(y)-R_{\rm in}] dy= 0 \qquad \forall x\in {\rm Supp}\; \overline n_\e, 
\eeq
\beq
\label{steady2}
\qquad a(x)+\f 1 \e \int K(x,y)[ \overline R _\e(y) -R_{\rm in}] dy \leq 0 \qquad \forall x\in \mathbb{R} \setminus {\rm Supp}\; \overline n_\e,
\eeq 
\beq \label{steady3}
\f{m(y)}{\e} [\overline R_{\rm in}(y) -  R_\e(y)] = \overline R_\e(y) \int K(x,y)\overline n_\e(x)dx \qquad \forall y \in \mathbb{R}.
\eeq
\end{definition}
We also recall the positivity condition which is important for the direct competition model \fer{comp_model} to have a unique ESD and to admit a Lyapunov functional \cite{PJ.GR:09}.
\begin{definition} [Operator positivity condition]
The function $c(x,y)$ satisfies a positivity condition, if 
 \beq \label{posit}\forall g\in M^1(\R)\setminus \{0\}, \qquad \int c(x,y)g(x)g(y)dxdy>0.\eeq
\end{definition}
For the model with the competition for resources a weaker assumption is enough to prove the uniqueness of ESD.
\begin{prop}[Uniqueness of ESD for \fer{full}] Assume \fer{as:ID}--\fer{as:coef2}.
Let $(\overline n_1,\overline R_1)$ and  $(\overline n_2,\overline R_2)$ be two evolutionary stable distributions for \fer{full}. Then for all $y\in \R$, $\overline R_1(y)=\overline R_2(y)$.
If moreover $K(x,y)$ is such that  
\beq \label{Knonz}\forall g\in M^1(\R)\setminus \{0\},\qquad \int K(x,y)g(x) dx \not \equiv 0,\eeq
then there exists at most one unique evolutionary stable distribution for \fer{full}.
\end{prop}
\proof
To prove the uniqueness of  ESD for \fer{full}, we follow the method used in \cite{PJ.GR:09,NC.PJ.GR:11} for the direct competition model.
Let $(\overline n_1,\overline R_1)$ and $(\overline n_2, \overline R_2)$ two evolutionary steady distributions. Using \fer{steady1}--\fer{steady2} we write
$$\int \overline n_1\left(a(x)+ \f 1 \e \int K(x,y) [\overline R_2(y)-R_{\rm in}] dy\right)dx\leq 0,$$
$$\int \overline n_2\left(a(x)+ \f 1 \e \int K(x,y) [\overline R_1(y)-R_{\rm in}] dy)\right)dx\leq 0,$$
$$-\int \overline n_i\left(a(x)+ \f 1 \e \int K(x,y) [\overline R_i(y)-R_{\rm in}] dy\right)dx= 0,\qquad \text{for }i=1,2.$$
By adding the above inequalities we obtain
$$\int \int  K(x,y)\left(\overline n_1(x)-\overline n_2(x)\right)\left(\overline R_2(y)-\overline R_1(y)\right)dydx\leq 0.$$
Replacing the values of $\overline R_1$ and $\overline R_2$ we obtain
$$\int \int  K(x,y)\f{m(y)}{\e}\left(\overline n_1(x)-\overline n_2(x)\right)\left(\f{R_{\rm in}(y)}{\f{m(y)}{\e}+\int K(x',y)\overline n_2(x')dx'}-\f{R_{\rm in}(y)}{\f{m(y)}{\e}+\int K(x',y)\overline n_1(x')dx'}\right)dydx\leq 0,$$
and thus
$$\int \f{\f{m(y)}{\e}R_{\rm in}(y)}{\left(\f{m(y)}{\e}+\int K(x',y)\overline n_2(x')dx'\right)\left(\f{m(y)}{\e}+\int K(x',y)\overline n_1(x')dx'\right)}\left(\int K(x,y)\left(\overline n_1(x)-\overline n_2(x)\right)dx\right)^2dy\leq 0.$$
It follows that
\beq\label{kn1}\int K(x,y) \overline n_1(x)dx = \int K(x,y) \overline n_2(x)dx,\qquad \text{a.e.}\eeq
and thus 
$$\overline R_1(y)=\overline  R_2(y), \qquad \text{a.e.}.$$
%
Moreover if $K(x,y)$ satisfies \fer{Knonz}, we obtain
$$  \overline n_1 - \overline n_2=0, \quad \text{in the sense of measures}.$$
\qed

For an ESD and following \cite{SM.BP.JW:10}, we define the Lyapunov functional for the model of competition for resources as
$$
S_{cr}(t)=-\int \overline n_\e(x) \ln n_\e(x,t)dx-\int \overline R_\e(y) \ln R_\e(y,t)dy+\int  n_\e(x,t)dx+\int   R_\e(y,t)dy.
$$
The dissipation associated with this functional is given by the inequality
\beq \begin{array}{rl}
  \f{dS_{cr}}{dt}(t):= D_{cr}(t)= &-\int \f{m(y) R_{\rm in}(y)}{\e^2 \overline R_\e(y) \; R_\e(y,t)}\left(\overline R_\e(y) - R_\e(y,t)\right)^2dy
  \\[4mm]
  & +\int n_\e(x,t)\left(  a(x) + \f 1 \e \int K(x,y)[\overline R_\e(y)- R_{\rm in}(y)]dy \right)dx \leq 0.
\end{array}
\label{eq:lyapunov}
\eeq
\\

Because we  expect that, according to the proof of Proposition \ref{prop:ntoc}, the expansion holds $R_\e (y,t)= R_{\rm in}(y) + O(\e)$, the limiting Lyapunov functional is reduced (up to an additive constant) to 
$$
S_{dc}(t)=-\int \overline n(x) \ln n(x,t)dx + \int n(x,t)dx. 
$$ 
This is indeed the Lyapunov functional used in for direct competition but $\f {d}{dt} S_{dc}(t)$ is negative only when the positivity condition \fer{posit} holds true. This condition appears naturally in our framework because we can compute the dissipation  associated with $
S_{dc}(t)$ 
$$
 \f{dS_{dc}}{dt}(t):= D_{dc}(t)
$$
A direct computation of $D_{dc}$ using \eqref{comp_model} is of course also possible but is less instructive than an expansion in $\e$ of $D_{cr}$ that uses simply 
$$
\overline R_\e(y) = R_{\rm in}(y) - \f{\e \overline R_\e(y) }{m(y)} \int K(x,y) \overline n_\e (x)dx.
$$
In the limit $\e \to 0$, we obtain
\beq \begin{array}{rl}
D_{dc}(t) =& -\int \f{R_{\rm in}(y)}{m(y)}  \left[ï¿½\int K(x,y) \big( n(x,t)- \overline n (x) \big)  \right]^2 dy
 \\[4mm]
  & + \int n(x,t)\left(a(x) - \int K(x,y) \f{R_{\rm in}(y)}{m(y)}\int K(x',y) \overline n(x')  dx'  dy \right)dx.
\end{array} 
\label{diis_dc} \eeq

Using the expression of $c(x,x')$ in \eqref{comp_kernel} yields the usual dissipation for the direct competition model 
$$\begin{array}{rl}
D_{dc} (t) = &- \int \int c(x,x') \big( n(x,t)- \overline n (x) \big) \big( n(x',t)- \overline n (x') \big)  dx\;  dx'
 \\[4mm]
& + \int n(x,t) \left( a(x) - \int   c(x,x') \overline n(x') dx' \right)dx.
\end{array} $$

The particular form of $c(\cdot,\cdot)$ in \eqref{comp_kernel} together with \fer{Knonz} ensures that the integral with the quadratic term in $n-\overline n$ is always positive (we referred to that as operator positivity) and the definition of an ESS makes that the linear term also gives a nonpositive contribution.

\section{Numerical illustration}
\label{sec:numerics}

How good is the  approximation of competition through resources by direct competition? Beyond the pure theoretical statement in section \ref{sec:ddc}, this can be illustrated thanks to some numerical simulations. We have performed such comparisons in  the case with mutations given by \eqref{full_mut}.  We assume that a resource uptake kernel and a distribution of resource supply are Gaussian.  A trait-dependent (and competition-independent) growth rate takes the maximum at $x=0$.
We have used the coefficients given as follows
\beq \label{data}
K(x,y)=\frac{1}{\sigma_K\sqrt{2\pi}}\exp\left(-\frac{(x-y)^2}{2\sigma_K^2}\right),\qquad R_{\rm in}(y)=\frac{M_{\rm in}}{\sigma_{\rm in}\sqrt{2\pi}}\exp(-\frac{y^2}{2\sigma_{\rm in}^2}), \qquad 
a(x)=1-x^2,
\eeq
and $m(y)$ is taken as a constant. See section \ref{sec:model} for the converted kernel of direct competition.  Throughout this section we have taken $\sigma_K=\sigma_{\rm in}=.5$, $\mu =0.005$ and we vary $\e$, $m$ and $M_{\rm in}$.
\\

\begin{figure}
{\centering
\includegraphics[width =8.5cm, height =5cm]{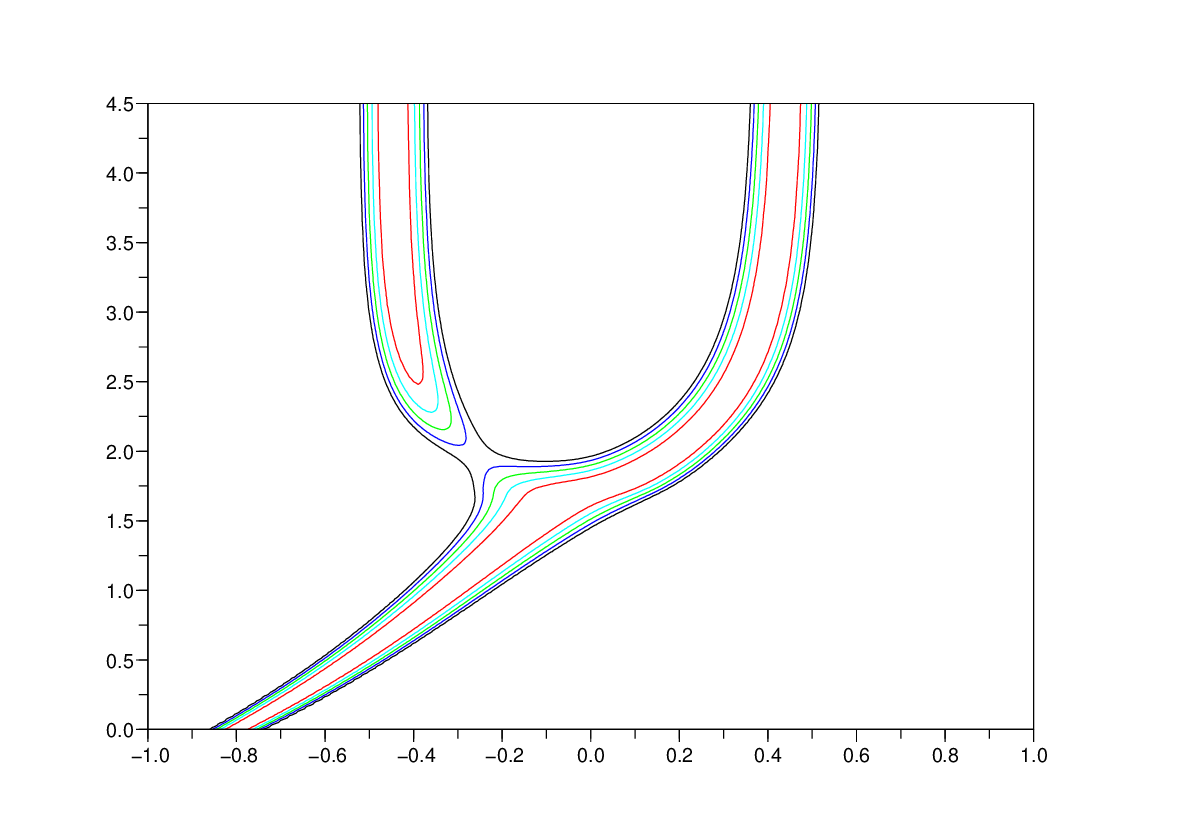}
\vspace{-.8cm} 
\caption{ \label{fig0} With our choice of parameters later (here $\e=0.001$, $m=10$, $M_{\rm in}=10$), the dynamics always undergo a branching from a monomorphic to a dimorphic population. This figure depicts the  typical contour lines of the solution to  \eqref {full_mut} (horizontal axis is $x$, vertical axis is time). 
} } \end{figure} 

The initial data $n^0$ is a Gaussian centered at $x=-0.8$ with its variance equal to the mutation rate $\mu$ (in accordance with the analysis in section \ref{mutations}) and $R$ is initially chosen equal to  $R_{\rm in}$. With all our choices of the remaining parameters $\e$, $m$, $M_{\rm in}$ below, the solution is always initially monomorphic and undergoes a dimorphic branching, then it stabilizes as depicted in Figure \ref{fig0}.
\\

\begin{figure}
{\centering
\includegraphics[width =8.5cm, height =6cm]{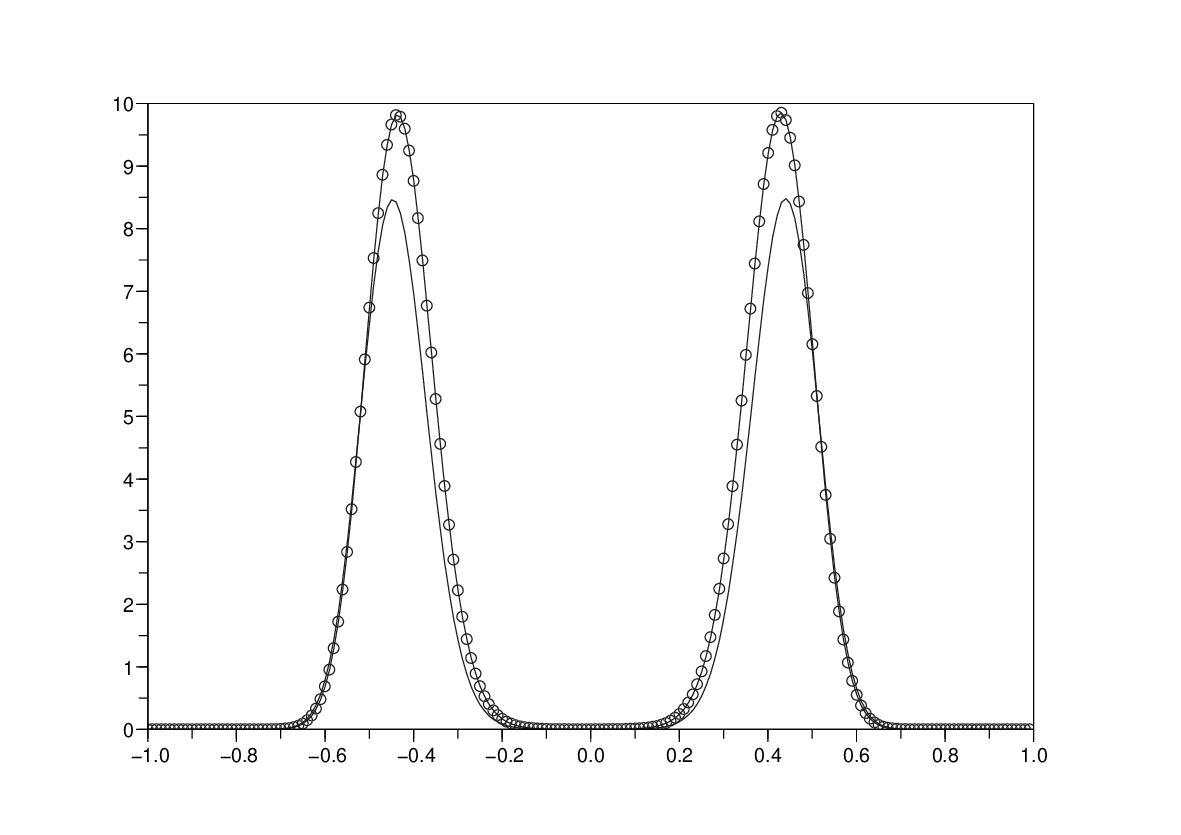} \hspace{-1cm}
\includegraphics[width =8.5cm, height =6cm]{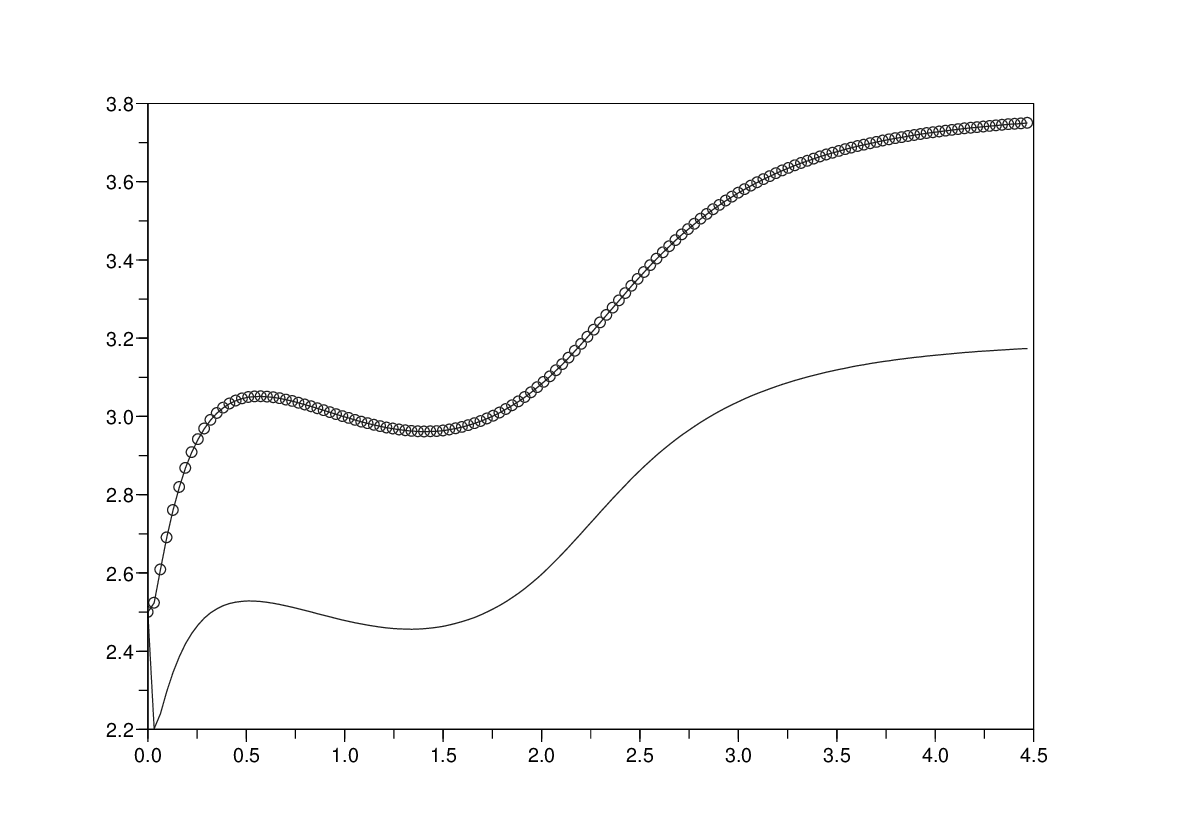}
\vspace{-.8cm} 
\caption{ \label{fig1} Comparisons between the model of competition for resources, with $m=M_{\rm in}=1$ and $\e=0.1$, (dotted curve) and the direct competition model (continuous line, in fact $\e=0.001$). Left: the dimorphic population density as a function of the trait $x$ for large times. Right: total population $\int n(x,t)dx$ as a function of time. 
} } \end{figure} 

We first compare the dynamics of the full model \eqref{full_mut} and that of the approximated model \eqref{comp_mut}.  To do so, we set $m= M_{\rm in}=1$ and studied two different values for the parameter $\e$, namely $\e=0.1$ and $\e=0.001$.  For $\e=0.1$, the relative error $\max \f{|R-R_{\rm in}|}{R_{\rm in}}=0.18$ was significantly large so the equilibrium distribution of resource  $R\approx R_{\rm in}$ was not achieved in this case.  As we decrease $\e$ to $\e=0.001$, we obtained $\max \f{|R-R_{\rm in}|}{R_{\rm in}}=0.0018$ and the resource distribution was well approximated by $R\approx R_{\rm in}$; thus, we expect that our reduction to direct competition model should give a good approximation.  Numerical comparisons for these cases are depicted in Figure \ref{fig1}, showing that the dimorphic distributions are well conserved in both solutions even though the total population is under-estimated in the direct competition model compared to $\e=0.1$ by a factor of approximately $25\%$.
\\

\begin{figure}
{\centering
\includegraphics[width =8.5cm, height =6cm]{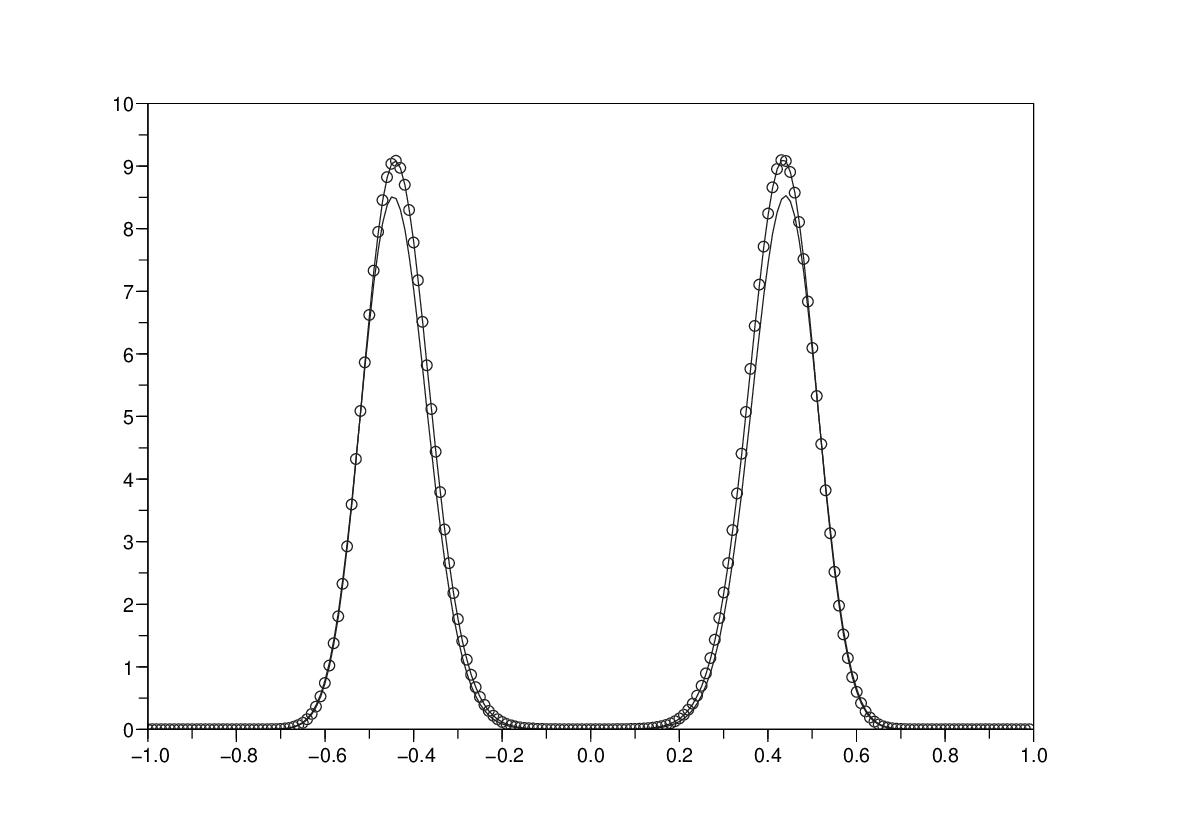} \hspace{-1cm}
\includegraphics[width =8.5cm, height =6cm]{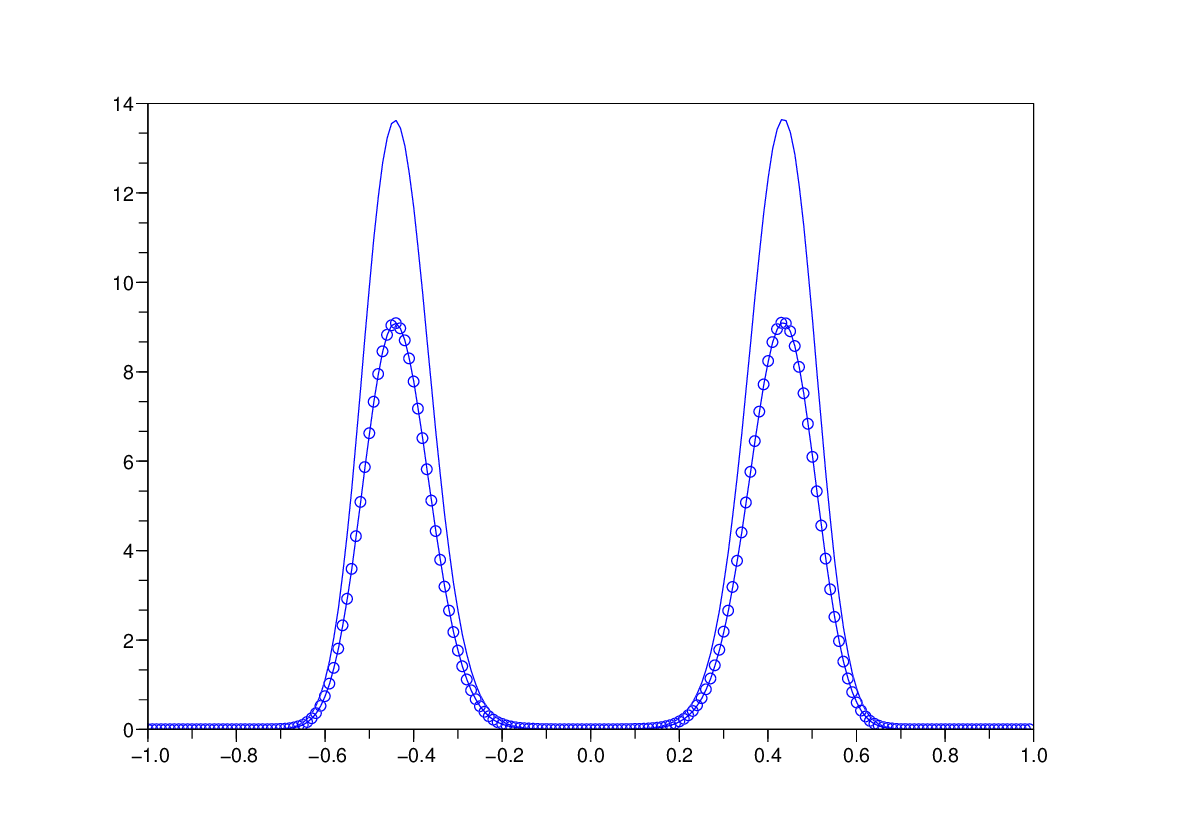}
\vspace{-.8cm} 
\caption{ \label{fig2} These figures illustrate that the direct competition kernel only depends on the ratio $R_{\rm in}/m$ and thus we can expect that the solutions to the model of competition for resource also does so.  (Left) we compare for $\e=0.05$, the population density solutions with $m=M_{\rm in}=1$  (dotted line) and with $m=M_{\rm in}=10$ (continuous line). (Right) we compare for $\e=0.05$ the solutions with 
$m=M_{\rm in}=1$  (dotted line) and with $m=1.5$, $M_{\rm in}=1$ (continuous line). When the ratio is preserved (even though coefficients differ by a factor $10$) the solutions are very close. But a small variation of the coefficients, that also changes the ratio, induces a visible difference on the solutions. 
} } \end{figure} 

The competition kernel \eqref{comp_kernel} in the direct competition approximation  has a remarkable property that it only depends upon the ratio $R_{\rm in}/m$. It implies that solution to \eqref{full_mut} depends only mildly on this ratio. This prediction from our analysis is confirmed by numerical calculations (Figure \ref{fig2}).

\section{The scales in McArthur's asymptotics}
\label{sec:macarthur}

After rescaling of \eqref{macarthur}, we consider the limit of 
\beq
\begin{cases}
\f{\p}{\p t}n_\e(x,t)=b(x) n_\e(x,t)\big[ \int K(x,y)R_\e (y,t) dy- d(x) \big],
\\[4mm]
\e  \f{\p}{\p t} R_\e (y,t)= R_\e(y,t) \big[m(y) (1 - \f{R_\e (y,t) }{R_{\rm cc}(y)}) - \int K(x,y)n_\e (x,t)dx \big],
\end{cases}
\label{ma}
\eeq

Putting $\e=0$,  McArthur in \cite{macarthur} builds on the formal limit $0=m(y) (1 - \f{R_\e (y,t) }{R_{cc}(y)}) - \int K(x,y)n_\e (x,t)dx $.  However, as noticed already in \cite{PA.CR:09,PA.CR.RD:08}, this is not correct because when this rate is negative, $R_\e(y,t)$ vanishes and thus the asymptotic expression is
$$
R (y,t) = R_{cc}(y) \left|  1-  \int \f{K(x,y)}{m(y)} n (x,t)dx  \right|_+ 
$$
where $| a |_+$ denotes the positive part of $a$. Inserting this formula in the equation for $n_\e$ gives us the asymptotic model
\beq
\f{\p}{\p t}n(x,t)=b(x) n(x,t)\big[ \int K(x,y) R_{cc}(y) \left|  1-  \int \f{K(x',y)}{m(y)} n (x',t)dx'  \right|_+ dy- d(x) \big].
\label{ma_comp}
\eeq
Only in the regime where resources are never exhausted (but there is no easy a priori characterization; this depends on the initial data and may change along the dynamics) one recovers direct competition \eqref{comp_model} with 
\beq
c(x, x') =  \int K(x,y)  \f{R_{\rm cc}(y)}{m(y)}  K(x',y) dy, \qquad a(x) = \f 1 \e  \int K(x,y) R_{cc}(y)  dy - d(x) .
\label{ma_cefp}
\eeq
As a consequence,  \cite{PA.CR:09,PA.CR.RD:08}  concluded that the behaviour of the resource-competition model and of the Lotka-Volterra one can differ significantly but did not attempt to clarify the conditions for the time-scaling argument to work. Here we propose another route which is to introduce another timescale for main resource consumption and guarantee it will not get extinct.
\bigskip

To fix this shortcoming, we may agin introduce two different timescales, as we proposed earlier. The rescaled version reads 
\beq
\begin{cases}
\f{\p}{\p t}n_\e(x,t)=b(x) n_\e(x,t)\big[ \f 1 \e \int K(x,y)R_\e (y,t) dy + a(x) - \f  1\e  \int K(x,y)R_{\rm cc} (y,t) dy \big],
\\[4mm]
\e^2  \f{\p}{\p t} R_\e (y,t)= R_\e(y,t) \big[m(y) (1 - \f{R_\e (y,t) }{R_{\rm cc}(y)}) - \e \int K(x,y)n_\e (x,t)dx \big],
\end{cases}
\label{maresc}
\eeq
Then the first order expansion in the equation on $R_\e$ is positive for $\e$ small enough 
$$
R_\e (y,t) = R_{cc}(y) \left[  1- \e \int \f{K(x,y)}{m(y)} n_\e (x,t)dx \right] + O(\e^2).
$$
Once inserted in equation \eqref{maresc}, we find the correct asymtotics
$$
\f{\p}{\p t}n_\e(x,t)=b(x) n_\e(x,t)\big[ a(x) -\int c(x, x')  n_\e(t,x')dx'  \big]
$$
with $c(x,x')$ as defined in \eqref{ma_cefp}.
\\

\noindent {\bf Acknowledgment} The three authors thank {\em Readilab}, LIA 197 CNRS,  for supporting their collaborations and also Meiji University and UPMC. J.Y.W. is partly supported by Glocal COE program "Formation and Development of Mathematical Sciences Based on Modeling and Analysis". S. M. benefits from a 2 year "Fondation Math\'ematique Jacques Hadamard" (FMJH) postdoc scholarship. She would like to thank Ecole Polytechnique for its hospitality.

{The authors would like also to thank G\'eza Mesz\'ena for careful review of this paper, valuable comments, several references and interpretations.  We borrowed some sentences  from the report. }

%
%


\end{document}